\documentclass[12pt]{amsart}

\usepackage{geometry}
\usepackage{amsmath,amssymb,amsthm}%
\usepackage{indentfirst,color,graphicx}%
\usepackage[colorlinks,citecolor=red,linkcolor=blue,urlcolor=cyan]{hyperref}%
\usepackage{tikz-cd}
\tikzcdset{arrow style=tikz,diagrams={>=stealth}}

\theoremstyle{plain}%
\newtheorem{theorem}{Theorem}[section]%
\newtheorem{corollary}[theorem]{Corollary}%
\theoremstyle{remark}%
\numberwithin{equation}{section}%

\newcommand{\BB}{\mathbb{B}}%
\newcommand{\CC}{\mathbb{C}}%
\newcommand{\RR}{\mathbb{R}}%
\newcommand{\calO}{\mathcal{O}}%
\newcommand{\pd}{\partial}%
\renewcommand{\geq}{\geqslant}%
\renewcommand{\leq}{\leqslant}%
\renewcommand{\Re}{\operatorname{Re}}%
\renewcommand{\tilde}[1]{\widetilde{#1}}%

\begin{document}

\title[A Converse of Berndtsson's Theorem]{A converse of Berndtsson's theorem on the positivity of direct images}

\author{Wang Xu}
\address{School of Mathematics, Sun Yat-sen University, Guangzhou 510275, China}
\email{xuwang@amss.ac.cn; xuwang6@mail.sysu.edu.cn}

\author{Hui Yang}
\address{School of Mathematical Sciences, Peking University, Beijing 100871, China}
\email{yanghui@amss.ac.cn}

\thanks{The first author is supported by National Key R\&D Program of China (No. 2024YFA1015200) and National Natural Science Foundation of China (No. 12501107).}

\maketitle

\begin{abstract}
Berndtsson's famous theorem asserts that, for a compact K\"ahler fibration $p:X\to Y$, the direct image bundle $p_*(K_{X/Y}\otimes L)$ of a semi-positive Hermitian holomorphic line  bundle $L\to X$ is Nakano semi-positive. As a continuation of our previous work, we prove a converse of Berndtsson's theorem in the case of a projective fibration: if $p_*(K_{X/Y}\otimes L\otimes E)$ is Griffiths semi-positive for every semi-positive Hermitian holomorphic line bundle $E\to X$, then the curvature of $L$ must be semi-positive.
\end{abstract}


\section{Introduction}

Positivity concepts, such as plurisubharmonicity, Griffiths positivity, and Nakano positivity, play an important role in several complex variables. In the $L^2$ theory for the $\bar{\pd}$-operator, numerous fundamental results are tied to these notions of positivity. For instance, given a plurisubharmonic function on a pseudoconvex domain, there are H\"ormander's $L^2$ existence theorem, Skoda's $L^2$ division theorem, the Ohsawa-Takegoshi $L^2$ extension theorem, as well as the optimal $L^2$ extension theorems established by B{\l}ocki and Guan-Zhou. These results serve as powerful tools in complex analysis, complex geometry, algebraic geometry, and other related fields.

Positivity is also the central theme in Berndtsson's ``\textit{complex Brunn-Minkowski theory}". In \cite{Berndtsson98} and \cite{Berndtsson06}, Berndtsson established complex analogues of the classical Pr\'ekopa theorem from convex analysis. In particular, he proved the plurisubharmonic variation of fibrewise Bergman kernels. From a higher viewpoint, Berndtsson \cite{Berndtsson09} further proved that the direct image bundle of a positive Hermitian holomorphic line bundle under certain holomorphic fibrations is Nakano positive. These contributions form a complex counterpart to the classical Brunn-Minkowski theory, with wide-ranging applications in complex analysis, K\"ahler geometry and algebraic geometry.

Berndtsson's theory is deeply connected to the $L^2$ theory. A crucial ingredient in his proofs is H\"ormander's $L^2$ estimate for $\bar{\pd}$. Moreover, Guan-Zhou \cite{GZ15} showed that their optimal $L^2$ extension theorem implies the plurisubharmonic variation of fiberwise Bergman kernels and the Griffiths positivity of direct images. Conversely, Berndtsson-Lempert \cite{BL16} derived an optimal $L^2$ extension theorem using the positivity of direct images.

In recent years, there are substantial research around the ``\textit{converse $L^2$ theory}" (see \cite[etc]{DWZZ, HI21, DNW, DNWZ}). This framework aims to characterize positivity notions via various $L^2$ conditions for $\bar{\pd}$ and thereby establish their necessity in $L^2$ theory. For example, Deng-Ning-Wang \cite{DNW} proved that an upper semi-continuous function satisfying the ``optimal $L^2$ extension condition" must be plurisubharmonic. More strikingly, Deng-Ning-Wang-Zhou \cite{DNWZ} showed the equivalence between Nakano positivity (a geometric property) and the ``optimal $L^2$ estimate condition" (an analytic property). This characterization also provides a new approach to the Nakano positivity of direct images (see \cite{DNWZ}).

Thus, positivity serves as the cornerstone uniting three closely interrelated theories: the $L^2$ theory, the converse $L^2$ theory, and the complex Brunn-Minkowski theory. In particular, we have the following diagram illustrating the situation:

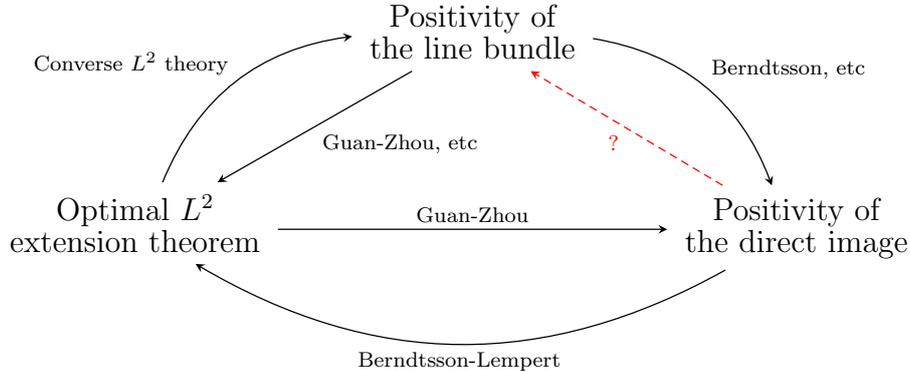
\begin{figure}[h]
\centering
\begin{tikzcd}
& \text{\normalsize Positivity of}\atop\text{\normalsize the line bundle} \arrow[ldd, "\text{Guan-Zhou, etc}"] \arrow[rdd, bend left, "\text{Berndtsson, etc}"] & \\
& & \\
\text{\normalsize Optimal $L^2$}\atop\text{\normalsize extension theorem} \arrow[rr, "\text{Guan-Zhou}"] \arrow[ruu, bend left, "\text{Converse $L^2$ theory}"] & & \text{\normalsize Positivity of}\atop\text{\normalsize the direct image} \arrow[ll, bend left, "\text{Berndtsson-Lempert}"] \arrow[luu, densely dashed, color=red, "?"]
\end{tikzcd}
\end{figure}

In view of this diagram, a natural question arises: does the positivity of the direct image imply the positivity of the original line bundle \textit{in some sense}? This question can be regarded as a parallel to the converse $L^2$ theory. However, only limited results have been obtained in this direction.\\

For clarity, we first recall the statement of Berndtsson's theorem.

\begin{theorem}[Berndtsson \cite{Berndtsson09}] \label{Thm:PosBundle}
	Let $p:X\to Y$ be a proper holomorphic submersion from a K\"ahler manifold $X$ to a connected complex manifold $Y$. Let $(L,h)$ be a semi-positive Hermitian holomorphic line bundle over $X$. Then the induced $L^2$ metric on the direct image bundle $p_*(K_{X/Y}\otimes L)$ is Nakano semi-positive.
\end{theorem}

In \cite{LiXuZhou}, the authors obtained a partial converse to the above theorem, using the method of Berndtsson-Lempert and the characterization of plurisubharmonic functions via the ``optimal $L^2$ extension condition".

\begin{theorem}[Li-Xu-Zhou \cite{LiXuZhou}] \label{Thm:LXZ}
Let $(L,h_L)$ be a Hermitian holomorphic line bundle over a projective manifold $\Omega$ with $H^0(\Omega,K_\Omega\otimes L)\neq\{0\}$. Let $m$ be a positive integer. Suppose that for any semi-positive Hermitian holomorphic line bundle $(E,h_E)$ on $X:=\mathbb{B}^m\times\Omega$, the direct image of $K_{X/\mathbb{B}^m}\otimes p_2^*L\otimes E$ over $\mathbb{B}^m$ is Griffiths semi-positive, where $p_2:X\to\Omega$ is the natural projection. Then the curvature of $(L,h_L)$ is semi-positive.
\end{theorem}

Therefore, the Griffiths positivity (a condition weaker than Nakano positivity) of a family of associated direct images implies the positivity of the original line bundle. By Berndtsson's theorem, it follows that these direct images are in fact Nakano semi-positive. However, Theorem \ref{Thm:LXZ} is ``\textit{extrinsic}" in nature, as it invokes line bundles on the external space $\mathbb{B}^m\times\Omega$. Clearly, an ``\textit{intrinsic}" version of the result is more desirable.

Note that the plurisubharmonic variation of fiberwise Bergman kernels is a specific instance of positivity for direct images. In our previous work \cite{XuYang}, using an approach different from \cite{LiZhou,LiXuZhou}, we obtained a converse result to Berndtsson's theorem concerning such variation.

\begin{theorem}[Xu-Yang \cite{XuYang}] \label{Thm:ConvBerndt}
Let $\Omega$ be a domain in $\CC_\tau^{m}\times\CC_z^n$ and $\varphi$ a strongly upper semi-continuous function on $\Omega$ such that $B_{\Omega_\tau}(z;e^{-\varphi_\tau})$ is not identically zero on $\Omega$. Assume that for any plurisubharmonic function $\psi\geq0$ on $\Omega$, the function
$$ (\tau,z) \mapsto \log B_{\Omega_\tau}(z;e^{-\varphi_\tau-\psi_\tau}) $$
is also plurisubharmonic on $\Omega$. Then $\varphi$ must be a plurisubharmonic function. Here, $\Omega_\tau:=\{z\in\CC^n:(\tau,z)\in\Omega\}$, $\varphi_\tau:=\varphi(\tau,\cdot)$ and $\psi_\tau:=\psi(\tau,\cdot)$.
\end{theorem}

As the statement invokes only functions on $\Omega$, this converse result is intrinsic. Moreover, we \cite{XuYang} provided counterexamples  showing that the additional plurisubharmonic weight $\psi$ is necessary.

The main goal of this article is to prove the following converse to Theorem \ref{Thm:PosBundle}, whose formulation is parallel with Theorem \ref{Thm:ConvBerndt}.

\begin{theorem} \label{MainThm}
Let $p:X\to Y$ be a holomorphic submersion from a projective manifold $X$ to a connected complex manifold $Y$. Let $(L,h_L)$ be a Hermitian holomorphic line bundle over $X$ such that $H^0(X_y,K_{X_y}\otimes L|_{X_y})\neq\{0\}$ for some $y\in Y$, where $X_y:=p^{-1}(y)$. Assume that for any semi-positive Hermitian holomorphic line bundle $(E,h_E)$ on $X$, the direct image sheaf $p_*(K_{X/Y}\otimes L\otimes E)$ is Griffiths semi-positive. Then the curvature of $(L,h_L)$ is semi-positive.
\end{theorem}

The conditions of Theorem \ref{MainThm} can be slightly weaken, but we adopt this formulation for simplicity. The direct image sheaf $\mathcal{E}:=p_*(K_{X/Y}\otimes L\otimes E)$ can be regarded as a bundle of vector spaces over $Y$, whose fiber at $y\in Y$ is
$$ \mathcal{E}_y := H^0(X_y,K_{X_y}\otimes L|_{X_y}\otimes E|_{X_y}). $$
The Hermitian metrics on $L$ and $E$ induce an $L^2$ metric on $\mathcal{E}$, defined by
$$ \|u_y\|^2 := \int_{X_y} |u_y|^2, \quad \forall u_y\in \mathcal{E}_y. $$
When $(L,h_L)$ has semi-positive curvature, $(\mathcal{E},\|\cdot\|)$ admits a natural structure as Hermitian holomorphic vector bundle \cite{Berndtsson09}. In general, $\mathcal{E}$ need not be locally free, and we shall define the Griffiths positivity in a \textit{singular} sense: $\mathcal{E}$ is called Griffiths semi-positive, if $\log\|\xi\|$ is plurisubharmonic for any local holomorphic section $\xi$ of $\mathcal{E}^*$. In particular, if $\mathcal{E}$ is Griffiths semi-positive, then the \textit{relative Bergman kernel metric} on $K_{X/Y}\otimes L\otimes E$ has semi-positive curvature (see Section 2 for details).

The proof combines the methods from \cite{LiXuZhou} and \cite{XuYang}, together with the argument by contradiction. Due to the absence of global coordinates and trivialization, the main difficulty lies in ``localization". Similar to Theorem \ref{Thm:LXZ}, the assumption that $X$ is \textit{projective} ensures the existence of a positive Hermitian holomorphic line bundle $A\to X$ with $H^0(X,A)\neq\{0\}$, which allows us to construct a semi-positive Hermitian holomorphic line bundle $E\to X$ and a global holomorphic section $f_E\in H^0(X,E)$ with prescribed local behaviour.

Combining Theorem \ref{MainThm} with Berndtsson's theorem yields the following immediate corollary. In particular, if all direct images in a special family are Griffiths semi-positive, then they are indeed Nakano semi-positive.

\begin{corollary}
Let $p:X\to Y$ be a holomorphic submersion from a projective manifold $X$ to a connected complex manifold $Y$. Let $(L,h_L)$ be a Hermitian holomorphic line bundle over $X$ such that $H^0(X_y,K_{X_y}\otimes L|_{X_y})\neq\{0\}$ for some $y\in Y$. Then the following statements are equivalent:

\begin{enumerate}
    \item the curvature of $(L,h_L)$ is semi-positive;
    \item for any semi-positive Hermitian holomorphic line bundle $E\to X$, the direct image $p_*(K_{X/Y}\otimes L\otimes E)$ is Nakano semi-positive.
    \item for any semi-positive Hermitian holomorphic line bundle $E\to X$, the direct image $p_*(K_{X/Y}\otimes L\otimes E)$ is Griffiths semi-positive.
\end{enumerate}
\end{corollary}

The remainder of the article is structured as follows: Section 2 provides necessary preliminaries, and Section 3 presents the proof of Theorem \Ref{MainThm}.\\

\textbf{Acknowledgments.} The authors sincerely thank their PhD supervisor Prof. Xiangyu Zhou for his generous help over the years. The authors also want to thank Zhi Li, Zhuo Liu and Xujun Zhang for useful discussions on related topics.

\section{Preliminaries}

In this section, we recall some preparatory materials concerning positivity concepts and direct images.

\textbf{Positivity concepts for line bundles.} Let $(L,h_L)$ be a Hermitian holomorphic line bundle over a complex manifold $X$. Let $U\subset X$ be an open set on which $L$ is trivial, and let $f_L$ be a holomorphic frame of $L|_U$. We write
$$ |f_L|_{h_L}^2=e^{-\varphi_L}, $$
then the Chern curvature of $(L,h_L)$ is locally given by $i\pd\bar\pd\varphi_L\otimes\textup{Id}_L$. Therefore, the Chern curvature of $(L,h_L)$ is semi-positive (resp. positive) if and only if the local weights $\varphi_L$ are psh (resp. strictly psh). In this case, we simply say the line bundle $(L,h_L)$ is \textit{semi-positive} (resp. \textit{positive}).

In general, for a holomorphic line bundle equipped with a \textit{singular} Hermitian metric, the curvature current is semi-positive if and only if the local weights of the metric are psh functions. In particular, we are interested in the following example. Given a Hermitian holomorphic line bundle $(L,h_L)\to X$, we define its Bergman kernel by
$$ B_L(z) := \sup\left\{ f(z)\otimes\overline{f(z)}: f\in H^0(X,K_X\otimes L), \int_X|f|_{h_L}^2\leq1 \right\}, $$
where $K_X$ is the canonical line bundle of $X$. Let $\mathbf{e}$ be a local frame of $K_X\otimes L$. We write $B_L = b\mathbf{e}\otimes\overline{\mathbf{e}}$, and define a possibly singular Hermitian metric $|\cdot|$ on $K_X\otimes L$ by setting $|\mathbf{e}|^2=1/b$. This metric is called the \textit{Bergman kernel metric}, and it has semi-positive curvature.

\textbf{Positivity concepts for vector bundles.} Let $(E,h_E)$ be a Hermitian holomorphic vector bundle of rank $r$ over a complex manifold $X$ of dimension $n$. Let $(U,(z_1,\ldots,z_n))$ be a coordinate chart of $X$ such that $E|_U$ is trivial, and let $(e_1,\ldots,e_r)$ be a holomorphic frame of $E|_U$. Locally, the Chern curvature tensor of $(E,h_E)$ can be written as
\begin{equation*}
	\Theta(E,h_E) = \sum_{j,k,\alpha,\gamma} R_{j\bar{k}\alpha}^\gamma dz_j\wedge d\bar{z}_k \otimes e_\alpha^* \otimes e_\gamma.
\end{equation*}
We define
\begin{equation*}
    R_{j\bar{k}\alpha\bar{\beta}} := \sum_\gamma R_{j\bar{k}\alpha}^\gamma h_{\gamma\bar{\beta}}, \quad \text{where } h_{\gamma\bar{\beta}} := \langle e_\gamma,e_\beta \rangle_{h_E}.
\end{equation*}
The vector bundle $(E,h_E)$ is said to be \textit{Griffiths semi-positive} (resp. positive) if
\begin{equation*}
	\sum_{j,k,\alpha,\beta} R_{j\bar{k}\alpha\bar{\beta}} a_j\overline{a_k}b_\alpha\overline{b_\beta} \geqslant 0 \text{ (resp. $>0$)}, \quad
	\forall~
	a\in\CC^n\setminus\{0\}, b\in\CC^r\setminus\{0\}.
\end{equation*}
Moreover, $(E,h_E)$ is said to be \textit{Nakano semi-positive} (resp. positive) if
\begin{equation*}
	\sum_{j,k,\alpha,\beta} R_{j\bar{k}\alpha\bar{\beta}} u_{j\alpha} \overline{u_{k\beta}} \geqslant 0 \text{ (resp. $>0$)}, \quad \forall~ u\in\CC^{nr}\setminus\{0\}.
\end{equation*}

These two positivity concepts coincide when $n=1$ or $r=1$. In general, Nakano positivity is strictly stronger than Griffiths positivity. It is well-known that $(E,h_E)$ is Griffiths semi-positive if and only if $\log|u|_{h_E^*}$ is psh for any nonvanishing local holomorphic section $u$ of the dual bundle $E^*$. This characterization also serves as a definition of Griffiths positivity for singular Hermitian metrics (see \cite{BP08}).

\textbf{Direct images.} Let $p:X\to Y$ be a proper holomorphic submersion between K\"ahler manifolds, then all the fibers $X_y:=p^{-1}(y)$ are smooth and compact. Let $(L,h_L)$ be a Hermitian holomorphic line bundle over $X$. We consider a bundle of finite dimensional vector spaces,
\begin{equation} \label{Eq:BundleE}
	E := \coprod_{y\in Y} E_y, \quad E_y := H^0(X_y, K_{X_y}\otimes L|_{X_y}).
\end{equation}
We denote by
$$ K_{X/Y}:=K_X-p^*K_Y $$
the relative canonical line bundle of $p$, then there exists a natural isomorphism
$$ K_{X_y}\cong K_{X/Y}|_{X_y}. $$
Therefore, the bundle $E$ is precisely the \textit{direct image} $p_*(K_{X/Y}\otimes L)$.

Assume that $(L,h_L)$ is semi-positive. By a variant of the Ohsawa-Takegoshi $L^2$ extension theorem, any element $u_y\in E_y$ can be extended to a holomorphic section $u\in H^0(p^{-1}(U),K_{X/Y}\otimes L)$, where $U\subset Y$ is a coordinate ball centered at $y$ (see \cite{Berndtsson09}). Starting from a basis of $E_y$, we obtain a local holomorphic frame for $E$, thus $E$ has a natural structure as a holomorphic vector bundle. Moreover, when we endow $E$ with the natural $L^2$ metric
$$ \|u_y\|^2 := \int_{X_y}|u_y|_{h_L}^2, \quad \forall u_y\in E_y, $$
Berndtsson's \cite{Berndtsson09} theorem asserts that $(E,\|\cdot\|)$ is Nakano semi-positive.

By identifying $K_{X_y}\otimes L|_{X_y}$ with $(K_{X/Y}\otimes L)|_{X_y}$, the fiberwise Bergman kernel metrics on $K_{X_y}\otimes L|_{X_y}$ glue into a metric on $K_{X/Y}\otimes L$, which is called the \textit{relative Bergman kernel metric}.

\textbf{Positivity concepts for direct images.} Let $p:X\to Y$ and $(L,h_L)$ be as above. In general, the bundle $E$ defined by \eqref{Eq:BundleE} need not be locally trivial, so we define a positivity concept for $E$ (or the direct image sheaf $p_*(K_{X/Y}\otimes L)$) in a singular sense.

A local section of $E$ on $U\subset Y$:
\[ U\ni y \mapsto u_y\in H^0(X_y,K_{X_y}\otimes L|_{X_y}), \]
is called \textit{holomorphic} if it is holomorphic as a section of $K_{X/Y}\otimes L$ over $p^{-1}(U)$. Moreover, a local section of $E^*:=\amalg_{y\in Y}(E_y)^*$ on $U\subset Y$:
\[ U\ni y \mapsto \xi_y\in \big(H^0(X_y,K_{X_y}\otimes L|_{X_y})\big)^*, \]
is called \textit{holomorphic} if $y\mapsto \xi_y(u_y)$ is holomorphic for any local holomorphic section $u$ of $E$. Finally, the bundle $E$ (or the direct image sheaf $p_*(K_{X/Y}\otimes L)$) is said to be \textit{Griffiths semi-positive} if $y\mapsto\log\|\xi_y\|$ is psh for any local holomorphic section $\xi$ of $E^*$, where $\|\cdot\|$ is the metric dual to the one on $E$. When $E$ is indeed a holomorphic vector bundle, this definition coincides with the usual one.

As explained in \cite{BP08}, the Griffiths semi-positivity of $p_*(K_{X/Y}\otimes L)$ would imply that the relative Bergman kernel metric on $K_{X/Y}\otimes L$ has semi-positive curvature:

Let $(\mathcal{N},(t,z))$ be a coordinate chart of $X$ such that $p|_{\mathcal{N}}$ is the trivial fibration $(t,z)\mapsto t$. We further assume that $\mathcal{N}$ is biholomorphic to a bounded product domain $U\times D$ in $\CC_t^m\times\CC_z^n$ and that $L|_{\mathcal{N}}$ is trivial. Let $\mathbf{e}$ be a holomorphic frame of $L$ on $\mathcal{N}$, then $dz\otimes\mathbf{e}$ serves as a holomorphic frame of $K_{X/Y}\otimes L$ on $\mathcal{N}$. We denote by $\psi(t,z)$ the local weight of the relative Bergman kernel metric, i.e.
$$ e^{\psi(t,z)} := \sup\left\{ |v_t(z)|^2 : u_t\in E_t, ~ \|u_t\|^2\leq1, ~ u_t|_D = v_t dz\otimes\mathbf{e} \right\}. $$
Given a holomorphic map $h:U\to D$, we can define a local holomorphic section $\xi$ of $E^*$ by
\begin{equation*}
	\xi_t(u_t) := v_t(h(t)), \quad \text{where } u_t\in E_t \text{ with } u_t|_D = v_tdz\otimes\mathbf{e}.
\end{equation*}
It is clear that
$$ \|\xi_t\|^2 = e^{\psi(t,h(t))}. $$
Since $p_*(K_{X/Y}\otimes L)$ is Griffiths semi-positive, we know that $\log\|\xi_t\|^2 = \psi(t,h(t))$ is psh. As $h:U\to D$ is arbitrary, we conclude that $\psi(t,z)$ is psh, and hence the relative Bergman kernel metric has semi-positive curvature.

\section{The Proof of Theorem \ref{MainThm}}

In this section, we prove a stronger version of Theorem \ref{MainThm}.

\begin{theorem} \label{Thm:ConvPosBundle}
	Let $p:X\to Y$ be a proper holomorphic submersion from a K\"ahler manifold $X$ of dimension $m+n$ to a connected complex manifold $Y$ of dimension $m$. Let $(L,h_L)$ be a Hermitian holomorphic line bundle over $X$. Assume that
	\begin{itemize}
		\item[(1)] $H^0(X_y,K_{X_y}\otimes L|_{X_y})\neq\{0\}$ for some $y\in Y$, where $X_y:=p^{-1}(y)$;
		\item[(2)] there exists a positive Hermitian holomorphic line bundle $(A,h_A)$ on $X$ with $H^0(X,A)\neq\{0\}$;
		\item[(3)] for any semi-positive Hermitian holomorphic line bundle $(E,h_E)$ on $X$, the relative Bergman kernel metric on $K_{X/Y}\otimes L\otimes E$ either has semi-positive curvature or is identically $+\infty$.
	\end{itemize}
	Then the curvature of $(L,h_L)$ is semi-positive.
\end{theorem}

The second assumption serves as an analogue of the existence of \textit{strictly} plurisubharmonic functions on $\CC^n$. A similar hypothesis appears in Theorem 1.1 and 1.2 of \cite{DNWZ}. Roughly speaking, having a positive Hermitian holomorphic line bundle, one can construct a \textit{global} semi-positive Hermitian holomorphic line bundle with a prescribed local weight (see Proposition 2.1 of \cite{DNWZ} and Step 3 of the following proof for details).

When $X$ is projective, condition (2) is automatically satisfied. Conversely, by the Kodaira embedding theorem, if $X$ is a compact manifold satisfying condition (2), then $X$ must be projective. As explained in Section 2, when the direct image $p_*(K_{X/Y}\otimes L\otimes E)$ is Griffiths semi-positive, then the relative Bergman kernel metric on $K_{X/Y}\otimes L\otimes E$ has semi-positive curvature. Therefore, Theorem \ref{MainThm} is a special case of Theorem \ref{Thm:ConvPosBundle}.\\

Notice that, since $X_y$ is compact and $h_L$ is smooth, we have
\[ \int_{X_y}|g|_{h_L}^2<+\infty, \quad \forall g\in H^0(X_y,K_{X_y}\otimes L|_{X_y}). \]
We denote by $h_B$ the relative Bergman kernel metric on $K_{X/Y}\otimes L$. By definition, $h_B(x)<+\infty$ if and only if there exists a holomorphic section 
$$ g\in H^0(X_y,K_{X_y}\otimes L|_{X_y}) \quad\text{with}\quad g(x)\neq0, $$
where $y:=p(x)$. Therefore, condition (1) implies $h_B\not\equiv+\infty$. Taking $(E,h_E)$ a trivial line bundle, it follows from condition (3) that the curvature of $h_B$ is semi-positive. Consequently, the local weights of $h_B$ are psh functions, and thus $\{h_B=+\infty\}$ is a pluripolar subset of $X$. We fix a positive Hermitian holomorphic line bundle $(A,h_A)$ and a nontrivial holomorphic section $f_A\in H^0(X,A)$ as assumed. Clearly,
$$ S:=\big\{f_A=0\big\}\bigcup\big\{h_B=+\infty\big\} $$
is a set of zero measure in $X$.

In the following, we prove the theorem by contradiction. For this purpose, we assume that
\begin{center}
	($\bigstar$) \quad the curvature of $(L,h_L)$ is NOT semi-positive at some point $x_0\in\Omega$.
\end{center}
Since $h_L$ is smooth, we may further assume that $x_0\notin S$.

For clearness, we divide the proof into four steps.\\

\textbf{\itshape Step 1: suitable trivialization around $x_0$.}

Since $p:X\to Y$ is a submersion, we can find a coordinate chart $(\mathcal{N},(\tau,z))$ around $x_0$ such that $\mathcal{N}$ is biholomorphic to a bounded product domain $U\times D$ in $\CC_\tau^m\times\CC_z^n$ and the restriction of $p$ to $\mathcal{N}$ is a trivial fibration: $(\tau,z)\mapsto\tau$. For simplicity, we always identify a point in $U\times D$ with the corresponding point in $\mathcal{N}\subset X$. Moreover, we consider $U$ as a subset of $Y$, and $D$ as a subset of $X_\tau$ for any $\tau\in U$. We will take
$$ dz := dz_1\wedge\cdots\wedge dz_n $$
as a holomorphic frame of $K_{X/Y}$ on $\mathcal{N}$.

Shrinking $\mathcal{N}$ if necessary, we may assume that $U\Subset Y$, $f_A\neq0$ on $U\times D$ and $L|_{U\times D}$ is trivial. By translation and scaling, we further assume that $x_0=(0,0)$ and $\BB^{2m}\times\BB^{2n} \Subset U\times D$. As $U\Subset Y$ and $p:X\to Y$ is proper, multiplying $f_A$ by a constant, we may assume that $|f_A|_{h_A}^2<1/e$ on $p^{-1}(U)$. Let
$$ e^{-\varphi_A} := |f_A|_{h_A}^2, $$
then $\varphi_A>1$ is a \textit{strictly} psh function on $U\times D$.

Let $f_L$ be a holomorphic frame of $L|_{U\times D}$ and we write
$$ e^{-\varphi_L} := |f_L|_{h_L}^2. $$
By assumption $(\bigstar)$, the weight function $\varphi_L$ is NOT psh at $(0,0)$. Then there exists some non-zero vector $(\eta,\xi)\in\BB^{2m}\times\BB^{2n}$ such that
$$ \phi(w) := \varphi_L(w\eta,w\xi) $$
is NOT subharmonic at $0\in\CC$. Since $\varphi_L$ is smooth, we may assume that $\eta\neq0$.\\

\textbf{\itshape Step 2: second-order estimate near $x_0$.}

In preparation for the forthcoming analysis, we establish some crucial estimates in this step.

We consider the Taylor expansion of $\phi(w)$ near $0\in\CC$:
$$ \phi(w) = \phi(0) + 2\Re\left( \frac{\pd\phi}{\pd w}(0)w + \frac{\pd^2\phi}{\pd w^2}(0)w^2 \right) - \alpha|w|^2+O(|w|^3), $$
where $\alpha:=-\tfrac{\pd^2\phi}{\pd w\pd\bar{w}}(0)$ is a \textit{positive} constant. Since
$$ \Re\left(\frac{\pd\phi}{\pd w}(0)w + \frac{\pd^2\phi}{\pd w^2}(0)w^2\right) $$
is a harmonic function, for any $r\in(0,1)$, we have
\begin{equation}\label{Eq:Est1}
	\begin{aligned}
		& \frac{1}{\pi r^2} \int_{\{w:|w|<r\}} \phi(w) d\lambda_w \\
		= &\, \frac{1}{\pi r^2} \int_{\{w:|w|<r\}} \left(\phi(0)-\alpha|w|^2+O(|w|^3)\right) d\lambda_w \\
		= &\, \phi(0) - \frac{\alpha}{2}r^2 + O(r^3).
	\end{aligned}
\end{equation}
We assume that the remaining term is bounded by $C_1r^3$ for $r\in(0,1)$, where $C_1$ is a positive constant. In particular, the mean-value inequality fails on sufficiently small disc.

As $\eta\neq0$, we may take a $\CC$-linear map $a:\CC^m\to\CC^n$ such that $a(\eta)=\xi$. Clearly, there exists some constant $R\in(0,1)$ such that
$$ \big\{ (\tau,z): |\tau|\leq R, |z-a(\tau)|\leq R \big\} \subset \BB^{2m}\times\BB^{2n}. $$
For any $\tau\in\BB^{2m}(0;R)$, we consider the Taylor expansion of $\varphi_L(\tau,\cdot)$ near $a(\tau)$:
\begin{multline*}
	\varphi_L(\tau,z) = \varphi_L(\tau,a(\tau)) + 2\Re\bigg( \sum_j \frac{\pd\varphi_L}{\pd z_j}(\tau,a(\tau))(z_j-a_j(\tau))\bigg) + O(|z-a(\tau)|^2).
\end{multline*}
Since $\varphi_L(\tau,z)$ is smooth, we may assume that the remaining term is bounded by $C_2|z-a(\tau)|^2$ on $\BB^{2n}(a(\tau);R)$ for any $\tau\in\BB^{2m}(0;R)$, where $C_2$ is a positive constant. For convenience, we denote
$$ \beta_\tau(z) := \sum_j \frac{\pd\varphi_L}{\pd z_j}(\tau,a(\tau))(z_j-a_j(\tau)), $$
which is a holomorphic function of $z$ vanishing at $a(\tau)$. In particular,
\begin{equation}\label{Eq:Est2}
	\big|\varphi_L(\tau,z) - \varphi_L(\tau,a(\tau)) - 2\Re\beta_\tau(z)\big| \leq C_2|z-a(\tau)|^2
\end{equation}
on $\BB^{2n}(a(\tau);R)$ for any $\tau\in\BB^{2m}(0;R)$.

Now, given $\tau\in\BB^{2m}(0;R)$, $r\in(0,R)$ and $v\in\calO(\BB^{2n}(a(\tau);r))$, it is clear that
\begin{equation}\label{Eq:Est3}
	\begin{aligned}
		&\, \int_{\{z:|z-a(\tau)|<r\}} |v(z)|^2e^{-\varphi_L(\tau,z)} d\lambda_z \\
		\geq &\, e^{-\varphi_L(\tau,a(\tau))-C_2r^2} \int_{\{z:|z-a(\tau)|<r\}} \big|v(z)e^{-\beta_\tau(z)}\big|^2 d\lambda_z \\
		\geq &\, \sigma_{2n}r^{2n} |v(a(\tau))|^2 e^{-\varphi_L(\tau,a(\tau))-C_2r^2},
	\end{aligned}
\end{equation}
where the last inequality follows from the plurisubharmonicity of $|ve^{-\beta_\tau}|^2$.

Recall that $x_0=(0,0)\notin S$, and then $h_B(x_0)<+\infty$. By the definition of $h_B$, there exists a holomorphic section
$$ g\in H^0(X_0,K_{X_0}\otimes L|_{X_0}) \quad\text{with}\quad g(0)\neq0, $$
where we identify $p(x_0)\in Y$ with $0\in U$. We write
$$ g|_D = udz\otimes f_L, $$
then $u$ is a holomorphic function on $D$ with $u(0)\neq0$. For later use, we consider the Taylor expansion of $u(z)e^{-\beta_0(z)}$ near $0\in\CC^n$:
\begin{equation*}
	u(z)e^{-\beta_0(z)} = u(0) + \sum_j c_jz_j + O(|z|^2),
\end{equation*}
where $c_j := \frac{\pd}{\pd z_j}\big|_{z=0}(ue^{-\beta_0})$. Consequently, for $r\in(0,1)$, we have
\begin{equation}\label{Eq:Est4}
	\begin{aligned}
		&\, \frac{1}{\sigma_{2n}r^{2n}} \int_{\{z:|z|<r\}} \big|u(z)e^{-\beta_0(z)}\big|^2 d\lambda_z \\
		= &\, \frac{1}{\sigma_{2n}r^{2n}} \int_{\{z:|z|<r\}} \bigg( |u(0)|^2 + \sum_j \overline{u(0)}c_jz_j + \sum_j u(0)\overline{c_jz_j} + O(|z|^2) \bigg) d\lambda_z \\
		= &\, |u(0)|^2 + O(r^2).
	\end{aligned}
\end{equation}
We assume that the remaining term is bounded by $C_3r^2$ for $r\in(0,1)$, where $C_3$ is a positive constant.\\

\textbf{\itshape Step 3: the construction of $(E,h_E)$ and $f_E\in H^0(X,E)$.}

In this step, we want to construct a semi-positive Hermitian holomorphic line bundle $(E,h_E)$ on $X$ together with a holomorphic section $f_E\in H^0(X,E)$ such that
\begin{gather*}
	|f_E|_{h_E}=1 \text{ on } \big\{ (\tau,z)\in U\times D: |\tau|\leq r, |z-a(\tau)|\leq\delta r \big\}, \\
	|f_E|_{h_E}<1 \text{ on } X_\tau\setminus\big\{ z\in D: |z-a(\tau)|\leq\delta r \big\} \text{ for any } \tau\in\overline{\BB^{2m}(0;r)},
\end{gather*}
where $r$ and $\delta$ are constants to be specified later.

Since we are dealing with smooth Hermitian metrics, we need a regularized max function. Let $\vartheta\in C_c^\infty(\RR)$ be a non-negative function with support in $[-\frac{1}{2},\frac{1}{2}]$ such that
\[ \int_\RR\vartheta(s)ds=1 \quad\text{and}\quad \int_\RR s\vartheta(s)ds=0. \]
We define
$$ \tilde{\max}\{t_1,t_2\} := \int_{\RR^2} \max\{t_1+s_1,t_2+s_2\} \vartheta(s_1)\vartheta(s_2) ds_1ds_2, \quad t_1,t_2\in\RR. $$
Then $\tilde{\max}\{t_1,t_2\}$ is non-decreasing in both variables, smooth and convex on $\RR^2$; $\max\{t_1,t_2\} \leq \tilde{\max}\{t_1,t_2\} \leq \max\{t_1,t_2\}+\frac{1}{2}$; if $t_1>t_2+1$, then $\tilde{\max}\{t_1,t_2\}=t_1$. Given a constant $T\in\RR$, for any $t_1<T-1$, it is clear that
$$ \big\{ t_2\in\RR: \tilde{\max}\{t_1,t_2\}\leq T \big\} = (-\infty,T]. $$

We choose a sufficiently small constant $\delta\in(0,1)$ such that
$$ \left(2C_2+\frac{C_3}{|u(0)|^2}\right)\delta^2 < \frac{\alpha}{6}, $$
where the constants $C_j$, $u(0)$ and $\alpha$ were introduced in Step 2.
By the properties of $\tilde{\max}$, it is elementary to check that the function
\begin{equation}
	\psi(\tau,z) := \tilde{\max} \left\{\log\frac{|\tau|}{3}, \log\frac{|z-a(\tau)|}{\delta}\right\}
\end{equation}
is psh on $U\times D$, and smooth on $(U\times D)\setminus\{(0,0)\}$. We choose a cutoff function $\rho\in C_c^\infty(U\times D)$ such that $0\leq\rho\leq1$ everywhere and $\rho\equiv1$ in a neighborhood of $(0,0)$. Since $\varphi_A$ is strictly psh, there exists an integer $q\gg1$ such that
$$ q\varphi_A + 2\rho\psi $$
is psh on $U\times D$. We may further assume that
$$q > \sup_{U\times D}\psi+1.$$
Since $U\times D$ is a bounded domain in $\CC^m\times\CC^n$, the right hand side is finite.

Now we consider the psh function
\begin{equation}
	\Psi := \tilde{\max} \big\{ q\varphi_A + 2\rho\psi, \psi \big\}
\end{equation}
on $U\times D$. Near the boundary of $U\times D$, we have
$$ q\varphi_A+2\rho\psi = q\varphi_A > q > \psi+1, $$
and then $\Psi=q\varphi_A$. Notice that $\psi$ has a log-pole at the origin. Therefore, in a sufficiently small neighborhood of $(0,0)$, we have
$$ q\varphi_A+2\rho\psi = q\varphi_A+2\psi \ll \psi, $$
and then $\Psi=\psi$. We choose $r\in(0,\frac{R}{3})$ such that $\Psi=\psi$ in a neighborhood of
$$ \big\{ (\tau,z): |\tau|\leq 3r, |z-a(\tau)|\leq \delta r \big\}. $$
For later use, we further require that
$$ C_1r < \frac{\alpha}{6}. $$

Let $T:=\log r<0$. We take a convex non-decreasing function $\chi\in C^\infty(\RR)$ such that $\chi\equiv0$ on $(-\infty,T]$, $\chi>0$ on $(T,+\infty)$ and $\chi(t)\equiv t-T-1$ for $t\geq T+2$. Then $\chi(\Psi)$ is a non-negative smooth psh function on $U\times D$. Near the boundary of $U\times D$, we have
$$ \Psi=q\varphi_A>q\gg1> T, $$
and then $\chi(\Psi)=q\varphi_A-T-1$. Recall that $f_A\in H^0(X,A)$ is non-vanishing on $U\times D$ and $|f_A|_{h_A}^2=e^{-\varphi_A}$. By taking $f_A^{\otimes q}$ as a local frame on $U\times D$, we can define a smooth Hermitian metric $h_E$ on $E:=A^{\otimes q}$ by
\begin{equation}
	\begin{cases}
		h_E := h_A^{\otimes q}e^{T+1}, & X\setminus(U\times D) \\
		|f_A^{\otimes q}|_{h_E}^2 := e^{-\chi(\Psi)}, & U\times D
	\end{cases}.
\end{equation}
Since $(A,h_A)$ is semi-positive and $\chi(\Psi)$ is psh, the curvature of $(E,h_E)$ is also semi-positive.

We consider the holomorphic section $f_E := f_A^{\otimes q} \in H^0(X,E)$. For convenience, we define
\begin{align*}
    W_1 := \big\{(\tau,z)\in U\times D: |\tau|\leq r,|z-a(\tau)|\leq\delta r\big\}, \\
    W_2 := \big\{(\tau,z)\in U\times D: |\tau|\leq r, |z-a(\tau)|>\delta r\big\}.
\end{align*}
On $W_1$, we have $\log\frac{|\tau|}{3}<T-1$ and $\log\frac{|z-a(\tau)|}{\delta}\leq T$, and then
$$ \Psi=\psi\leq T \quad\Rightarrow\quad \chi(\Psi)=0 \quad\Rightarrow\quad |f_E|_{h_E}^2=1. $$
On $W_2$, we have $\log\frac{|z-a(\tau)|}{\delta}>T$, and then
$$ \Psi\geq\psi>T \quad\Rightarrow\quad \chi(\Psi)>0 \quad\Rightarrow\quad |f_E|_{h_E}^2<1. $$
Recall that $|f_A|_{h_A}^2<1/e$ on $p^{-1}(U)$. For any $\tau\in\overline{\BB^{2m}(0;r)}$, it is clear that
$$ |f_E|_{h_E}^2 = |f_A|_{h_A}^{2q}e^{T+1} \leq e^{-q+T+1} < 1 \quad\text{on}\quad X_\tau\setminus D. $$
In summary, $(E,h_E)$ and $f_E\in H^0(X,E)$ satisfy all the requirements.\\

\textbf{\itshape Step 4: reaching a contradiction.}

In this final step, we aim to derive a contradiction by using the estimates proved in Step 2 and the line bundle $(E,h_E)$ constructed in Step 3.

As $H^0(X,E)\neq\{0\}$, it follows from condition (1) that the relative Bergman kernel metric $B_k$ on $K_{X/Y}\otimes L\otimes E^{\otimes k}$ is not identically $+\infty$. Moreover, by condition (3), the curvature current of $B_k$ is semi-positive. Since $f_E\neq0$ on $U\times D$, we take $dz\otimes f_L\otimes f_E^{\otimes k}$ as a holomorphic frame of $K_{X/Y}\otimes L\otimes E^{\otimes k}$ on $U\times D$. We denote by $\Phi_k$ the local weight of $B_k$, i.e.
\begin{equation*}
	e^{\Phi_k(\tau,z)} := \sup\left\{ \frac{|v(z)|^2}{\|f\|^2}: \begin{gathered}
		f\in H^0(X_\tau, K_{X_\tau}\otimes(L\otimes E^{\otimes k})|_{X_\tau}) \\ f\not\equiv0, ~ f|_D=vdz\otimes f_L\otimes f_E^{\otimes k}
	\end{gathered} \right\},
\end{equation*}
then $\Phi_k$ is a psh function on $U\times D$.

Recall that $|f_L|_{h_L}^2 = e^{-\varphi_L}$ on $U\times D$ and
$$ |f_E|_{h_E}\equiv1 \quad\text{on}\quad \big\{ (\tau,z): |\tau|\leq r, |z-a(\tau)|\leq\delta r \big\}. $$
Given $\tau\in\BB^{2m}(0;r)$, for any holomorphic section
$$ f\in H^0(X_\tau, K_{X_\tau}\otimes(L\otimes E^{\otimes k})|_{X_\tau}) \quad\text{with}\quad f|_D=vdz\otimes f_L\otimes f_E^{\otimes k}, $$
it follows from the estimate \eqref{Eq:Est3} that
\begin{align*}
	\|f\|^2 & \geq \int_{\{z:|z-a(\tau)|<\delta r\}} |f|^2 \\
	& = \int_{\{z:|z-a(\tau)|<\delta r\}} |v(z)|^2e^{-\varphi_L(\tau,z)} d\lambda_z \\
	& \geq \sigma_{2n}(\delta r)^{2n} |v(a(\tau))|^2 e^{-\varphi_L(\tau,a(\tau))-C_2(\delta r)^2}.
\end{align*}
As a consequence, for any $\tau\in\BB^{2m}(0;r)$, we have
\begin{equation} \label{Eq:Upper}
	\Phi_k(\tau,a(\tau)) \leq \varphi_L(\tau,a(\tau)) + C_2(\delta r)^2 - \log \big(\sigma_{2n}(\delta r)^{2n}\big).
\end{equation}

Recall that, as $x_0\notin S$, there is a holomorphic section
$$ g\in H^0(X_0,K_{X_0}\otimes L|_{X_0}) \quad\text{with}\quad g|_D=udz\otimes f_L \quad\text{and}\quad u(0)\neq0. $$
We consider the holomorphic section
$$ f_k := g\otimes f_E^{\otimes k} \in H^0(X_0,K_{X_0}\otimes (L\otimes E^{\otimes k})|_{X_0}). $$
Since
\begin{equation*}
    \begin{cases}
        |f_E|_{h_E}\equiv1 \text{ on } \{z:|z|\leq\delta r\}\subset X_0, \\
        |f_E|_{h_E}<1 \text{ on } X_0\setminus\{z:|z|\leq\delta r\},
    \end{cases}
\end{equation*}
it follows from the Lebesgue dominated convergence theorem that
\begin{equation*}
	\lim_{k\to+\infty} \|f_k\|^2 = \int_{\{z:|z|\leq\delta r\}} |g|^2 = \int_{\{z:|z|<\delta r\}} |u(z)|^2e^{-\varphi_L(0,z)} d\lambda_z.
\end{equation*}
According to the estimates \eqref{Eq:Est2} and \eqref{Eq:Est4},
\begin{align*}
	&\, \int_{\{z:|z|<\delta r\}} |u(z)|^2e^{-\varphi_L(0,z)} d\lambda_z \\
	\leq &\, \int_{\{z:|z|<\delta r\}} |u(z)e^{-\beta_0(z)}|^2 e^{-\varphi_L(0,0)+C_2(\delta r)^2} d\lambda_z  \\
	\leq &\, \sigma_{2n}(\delta r)^{2n} \big(|u(0)|^2+C_3(\delta r)^2\big) e^{-\varphi_L(0,0)+C_2(\delta r)^2}.
\end{align*}
By the definition of $\Phi_k$,
\begin{align*}
	\varliminf_{k\to+\infty} e^{\Phi_k(0,0)} & \geq \lim_{k\to+\infty} \frac{|u(0)|^2}{\|f_k\|^2} \\
	& \geq \frac{1}{\sigma_{2n}(\delta r)^{2n} e^{-\varphi_L(0,0)+C_2(\delta r)^2}} \frac{|u(0)|^2}{|u(0)|^2+C_3(\delta r)^2},
\end{align*}
and then
\begin{equation} \label{Eq:Lower}
	\varliminf_{k\to+\infty} \Phi_k(0,0) \geq \varphi_L(0,0) - C_2(\delta r)^2 - \frac{C_3}{|u(0)|^2} (\delta r)^2 - \log\big(\sigma_{2n}(\delta r)^{2n}\big).
\end{equation}

Since $\Phi_k$ are psh functions on $U\times D$, we have the mean-value inequality
$$ \Phi_k(0,0) \leq \frac{1}{\pi r^2} \int_{\{w:|w|<r\}} \Phi_k(w\eta,w\xi) d\lambda_w. $$
Combining the inequalities \eqref{Eq:Upper} and \eqref{Eq:Lower}, we have
\begin{align*}
	&\, \varphi_L(0,0) - C_2(\delta r)^2 - \frac{C_3}{|u(0)|^2}(\delta r)^2 - \log\big(\sigma_{2n}(\delta r)^{2n}\big) \\
	\leq &\, \frac{1}{\pi r^2} \int_{\{w:|w|<r\}} \left( \varphi_L(w\eta,w\xi) + C_2(\delta r)^2 - \log \big(\sigma_{2n}(\delta r)^{2n}\big) \right) d\lambda_w.
\end{align*}
Equivalently,
\begin{equation*}
	\phi(0) \leq \frac{1}{\pi r^2} \int_{\{w:|w|<r\}} \phi(w)d\lambda_w + \left(2C_2+\frac{C_3}{|u(0)|^2}\right) \delta^2r^2.
\end{equation*}
Combining the inequality \eqref{Eq:Est1}, we get
$$ \frac{\alpha}{2}r^2 \leq \left(2C_2+\frac{C_3}{|u(0)|^2}\right) \delta^2r^2 + C_1 r^3. $$
As $\delta$ and $r$ are chosen such that
\begin{equation*}
	\left(2C_2+\frac{C_3}{|u(0)|^2}\right)\delta^2 < \frac{\alpha}{6} \quad\text{and}\quad C_1r < \frac{\alpha}{6},
\end{equation*}
we reach a contradiction under the assumption $(\bigstar)$.

This completes the proof of Theorem \ref{Thm:ConvPosBundle}.

\end{document}